\documentclass[10pt,a4paper]{article}
\usepackage[T1]{fontenc}
\usepackage{amsmath}
\usepackage{amsfonts}
\usepackage{amssymb}
\usepackage{graphicx}
\usepackage{caption}
\usepackage{csquotes}
\usepackage[english]{babel}

\usepackage[backend=biber, style=alphabetic]{biblatex}
\graphicspath{{./images/}}

\title{Insolvability of $x^x = a$ in elementary functions\thanks{The paper was supported by Russian science foundation grant No 22-19-20073}}
\author{Alexey Kanel-Belov, Rodion Zaytsev}
\newcommand{\Addresses}{{
  \bigskip
  \footnotesize

  Rodion Zaytsev, \textsc{
Faculty of Mathematics, National Research University Higher School of Economics, Usacheva str. 6, Moscow, 119048, Russian Federation;\smallskip\newline 
Igor Krichever Center for Advanced Studies, Skolkovo Institute of Science and Technology, Bolshoy Boulevard 30, bld. 1, Moscow, 121205, Russian Federation}
\par\nopagebreak
 \textit{E-mail address}: \texttt{rvzaytsev@edu.hse.ru}

}}

\date{November 2023}

\begin{document}
\maketitle
\begin{abstract}
\noindent In paper \cite{kanel2019solvability} the insolvability in elementary functions of equation $\tan(x) - x = a$ was proved. This work applies the same topological method to prove the insolvability of equation $x^x = a$.
\end{abstract}
\section{Introduction}

In \cite{kanel2019solvability} the technique from topological Galois theory informally discussed in the following paragraph was applied to prove the insolvability in elementary functions of the equation $\tan(x) - x = a$. 

Suppose we are trying to solve in elementary functions (i.e. roughly speaking compositions  and algebraic operations of the logarithm, exponent and identity functions) an equation $f(x) = a$, where $f(x)$ is holomorphic and non-constant. When $a$ draws a closed curve in the complex plane, the roots also draw curves which, however, aren't necessarily closed. A permutation on roots is thus induced, and it is shown in \cite{kanel2019solvability} that if this group is not solvable, then the equation cannot be solved in elementary functions.

\section{Preliminary observations}
In this section the equation is transformed to make it easier to analyse, and then the standard properties, such as the location of the critical points, are derived.
\subsection{Transform}
Taking logarithm twice on both sides, we obtain an equivalent equation
$$ \ln(x) + \ln\ln(x) = \ln\ln(a) $$
Since we only care about the solvability in elementary functions, any transforms of $x, a$ which can be expressed as compositions elementary functions are allowed (because we can express the initial parameters through the transformed in elementary functions). Let 
$$z = \ln\ln(x), b = \ln\ln(a)$$
Then we have an equivalent equation
$$z + e^z = b$$
To make the notation consistent with \cite{kanel2019solvability}, rename $ b\to a$, so that the equation is 
$$z + e^z = a$$
\subsection{Critical points}
Recall that for a complex function $f(z)$ multiple roots merging is equivalent to $f'(z) = 0$. In our case, 
$$
f(z) = z + e^z \Rightarrow f'(z) = 1 + e^z
$$
$$
f'(z) = 0 \Leftrightarrow e^z = -1 \Leftrightarrow z = (2n+1)\pi i, n \in \mathbb{Z}
$$
Notice that the critical points are all located along the imaginary axis.
The corresponding values of $a$ are 
$$
a = f((2n+1)\pi i) = (2n+1)\pi i - 1
$$
Thus the corresponding $a$-s are just shifted to the left. We shall call the $n$-th critical point $z_n$, and the corresponding $a(z_n) = a_n$ (where $n$ is the parameter that appears naturally above). Let's find the order of the critical points.
$$f''(z) = e^z \neq 0$$
So all the critical points are precisely of the first order, that is there are exactly two roots that are merged.
\section{Permutation}
In this section, the connection between the path and the permutation is investigated.
\subsection{Path choice}
First of all, for the sake of simplicity, let's consider $a = 0$ the starting point.
\begin{center}
	\includegraphics[width=\textwidth]{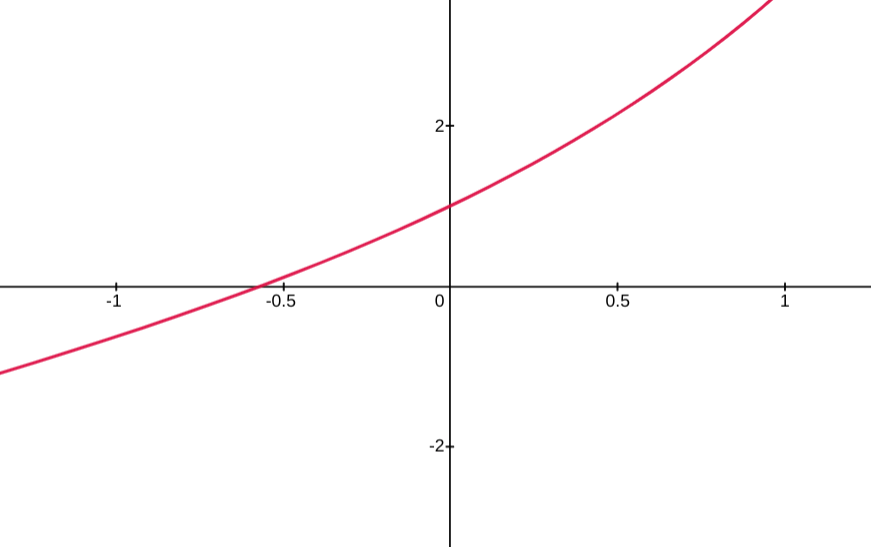}
	\captionof{figure}{$y = x+e^x$}
\end{center}
It can be easily shown that there is only one real root, as seen on the graph. Let's consider the path $z(t)$, which brings this real root to the $n$-th critical point, then makes a semi-loop around it, exchanging it with the other root, and then the new root comes back along the same path. For the sake of simplicity, case $n=2$ is illustrated.
\begin{center}
	\includegraphics[width=\textwidth]{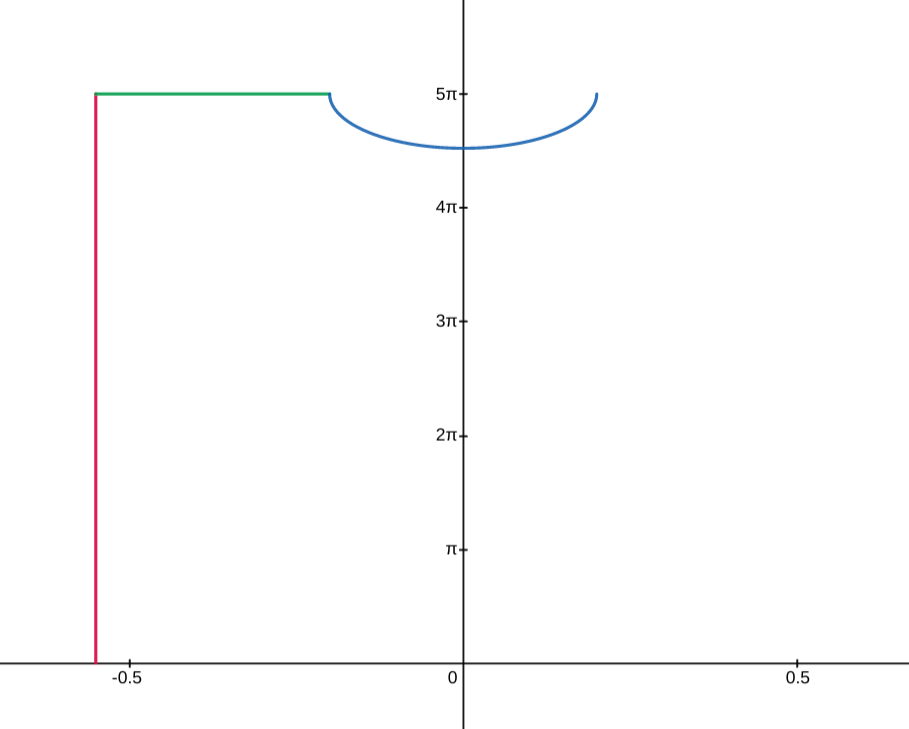}
	\captionof{figure}{The path of the real root}
\end{center}
\subsection{Path calculation}
We will now derive the path that $a$ follows as the real root takes the path described above. In the first part of the path, 
$$ z(t) = x + it$$
where $x$ is the real root (i.e. $e^{x} = -x$), so
$$ a(t) = x + it + e^{x + it} = x(1-\cos t) + i(t - x\sin t)$$ An important observation to make is that the cycles bend around the 'dangerous' points, that is $a_k$, because when $t = (2k+1)\pi$, 
$$a = 2x + i(2k+1)\pi$$
which is to the left of $a_k$, since $x < -\frac{1}{2}$ (this can be easily proved), and the real part of $a_k$ is $-1$. This also shows that by the end of the first part of the path, $a$ is located to the left of $a_n$. Let $y_n = (2n+1)\pi$, then the second part of the path is given by 
$$ z(s) = s+iy_n + e^{s+iy_n} \Rightarrow a(s) = s - e^s + iy_n$$
Thus $a(s)$ also goes in a straight line parallel to the real axis. As for the third part of the path, when the two roots are swapped $a$ makes a loop around $a_n$, as follows from the general theory. We can therefore sketch the path drawn by $a$
\begin{center}
	\includegraphics[width=\textwidth]{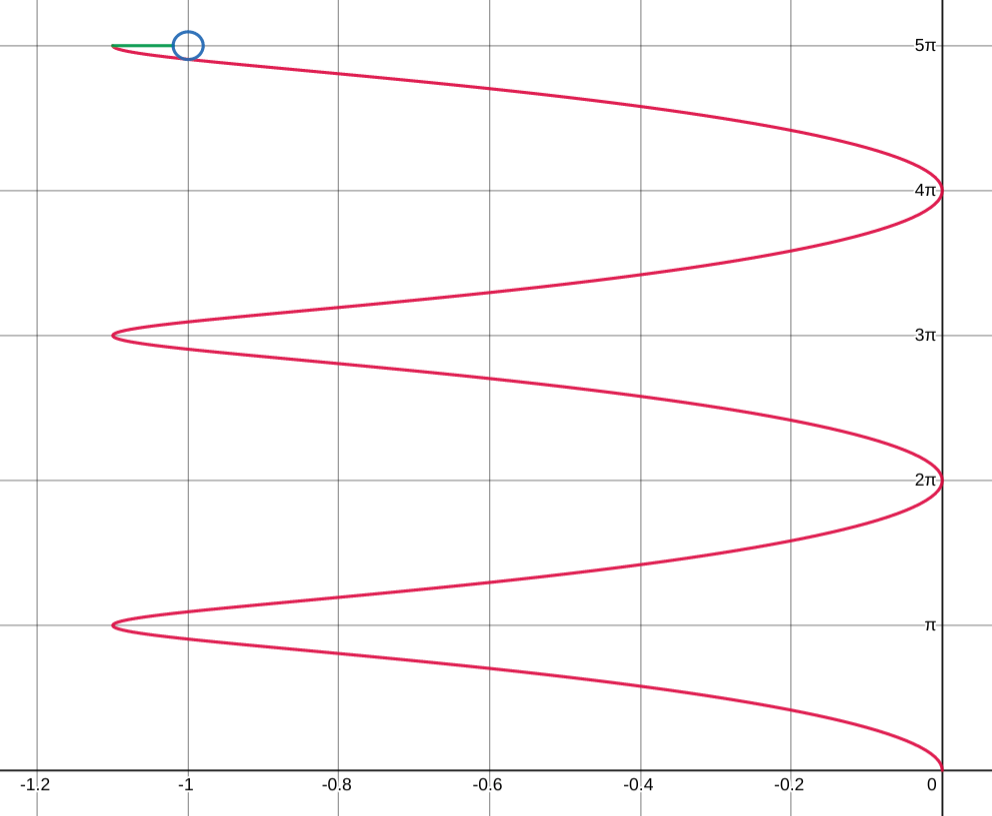}
	\captionof{figure}{The path followed by $a$}

\end{center}
This path is homotopically equivalent to a simple loop around the corresponding point.
\begin{center}
	\includegraphics[scale=0.3]{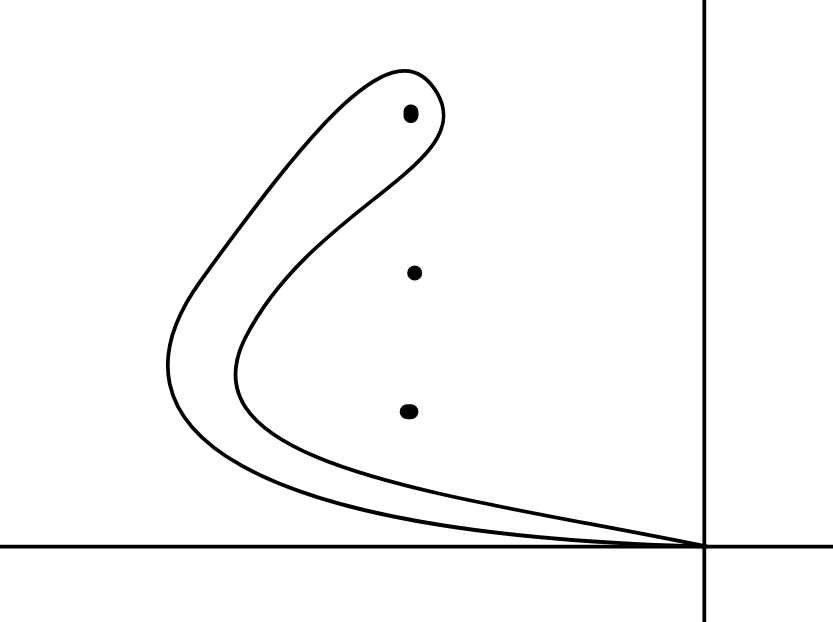}
	\captionof{figure}{Homotopically equivalent loop}

\end{center}
\pagebreak
\section{Main result}
In this section the insolvability of the given equation is proved.
Since a plane with a discrete set of points removed is homotopically equivalent to a bouquet of circles, it follows that any loop around the given points can be decomposed into a product of simple loops, where a simple loop looks like the one above.
This means that any possible permutation can be achieved by subsequently transposing the real root with one of the others multiple times. On the other hand, we know that the group is transitive, because it is true in the general case. Consider now the decomposition of permutation $\sigma$ that makes $r_n$ the real root (number the roots with natural numbers in some way such that $r_1 = x$ is the real root. As shown above, it will have the form 
$$ \sigma = (1n_1)(1n_2)\ldots(1n_N)$$
Clearly, the transposition $(1n)$ must appear in this product, otherwise $r_n$ will remain untouched. Hence the group contains all the transpositions $(1n)$ and is therefore the symmetric group on all the roots. In particular, the symmetric group is not solvable (as the number of roots is infinite). Therefore the equation is not solvable in elementary functions.
\printbibliography

\Addresses
\end{document}